\documentclass[12pt]{article}
\usepackage{amsmath}
\usepackage{amsthm}
\usepackage{amssymb}
\usepackage{mathrsfs}
\usepackage{authblk}
\usepackage[usenames]{color}
\usepackage[latin5]{inputenc}
\usepackage{cite}
\usepackage{hyperref}
\usepackage{pdfsync}
\usepackage{enumitem}
\setlist{parsep=0pt,itemindent=0pt}

\theoremstyle{definition}

\numberwithin{equation}{section}
\numberwithin{thm}{section}
\numberwithin{lemma}{section}
\numberwithin{prop}{section}
\numberwithin{cor}{section}
\numberwithin{rmk}{section}
\numberwithin{defn}{section}
\numberwithin{exa}{section}

\newcommand{\gen}[1]{\frac{\partial}{\partial{#1}}}
\newcommand{\pr}[1]{{\rm pr}^{(#1)}}
\newcommand{\curl}[1]{ \left\{#1\right\} }

\def\bu{\mathbf{u}}
\def\be{\mathbf{e}}
\def\<#1>{\langle#1\rangle}

\def\pd#1#2{\frac{\partial #1}{\partial #2}}

\DeclareMathOperator{\aff}{\mathfrak{aff}}

\DeclareMathOperator{\Sl}{\mathfrak{sl}}

\DeclareMathOperator{\SL}{SL}

\begin{document}
\pagenumbering{arabic}
\clearpage
\thispagestyle{empty}

\title{Invariance of second order ordinary differential equations under two-dimensional affine subalgebras of Ermakov--Pinney Lie algebra}

\author[1]{J.~F. Cari\~nena\thanks{jfc@unizar.es}}
\author[2]{F.~G\"ung\"or\thanks{gungorf@itu.edu.tr}}
\author[3]{P.~J. Torres\thanks{ptorres@ugr.es}}

\affil[1]{Departamento de F\'{\i}sica Te\'orica,
Universidad de Zaragoza, 50009 Zaragoza, Spain}
\affil[2]{Department of Mathematics, Faculty of Science and
Letters, Istanbul Technical University, 34469 Istanbul, Turkey}
\affil[3]{Departamento de Matem\'atica Aplicada, Universidad de
Granada, 18071 Granada, Spain}

\date{\today}

\maketitle




\begin{abstract}
Using the only admissible rank-two realisations of the Lie algebra of the affine group in one dimension in terms of the Lie algebra of Lie symmetries of the Ermakov-Pinney (EP) equation, some classes of second order nonlinear ordinary differential equations solvable by reduction method are constructed. One class includes the standard EP equation as a special case. A new EP equation with a perturbed potential but admitting the same solution formula as EP itself arises. The solution of the dissipative EP equation is also discussed.
\end{abstract}

\section{Introduction}

Among the most general second order linear differential equations in normal form
\begin{equation}
\psi''+q(x)\psi'+ p(x)\psi=0,
\label{eq_sode}
\end{equation}
those with $q\equiv0$ are especially important in both classical and quantum physics and will be said to be of
Schr\"odinger type, because the usual Schr\"odinger equation for the determination of stationary states is of this type
with a coefficient $p$ given by $p(x)=E-V(x)$, where $V$ is the potential and $E$ is the energy eigenvalue. Then, the equation can be written as
\begin{equation}
\psi''+ p(x) \psi=0\label{Hill}.
\end{equation}
Changing variables to $t$ for the independent variable and $x$ for the dependent variable, and using the dot notation for time derivative, the corresponding equation
\begin{equation}
\ddot x+p(t)x=0\label{oscil}
\end{equation}
is known in a more mathematical context as {\em Hill's equation} \cite{Hill,MW,GT12}. It has been shown in \cite{Belmonte-BeitiaPerez-GarciaVekslerchikTorres2007}
(see also \cite{Belmonte-BeitiaPerez-GarciaVekslerchikTorres2008}) that there is an infinitesimal point transformation  of symmetry of such a
Schr\"odinger type equation  of the form
\begin{equation}\label{EP-svfa}
  X_a(t,x)=a(t)\pd{} t+\frac{\dot{a}(t)}{2}x\pd{} x,
\end{equation}
where the function $a$ satisfies the following third order linear ODE
\begin{equation}\label{3rd2}
 \mathbb{M}(a)=\dddot{a}+4p(t)\dot{a}+2\dot{p}(t) a=0,
\end{equation}
which was called in \cite{Gu02} projective vector field equation. Moreover, as the differential equation (\ref{oscil}) is linear,  all vector fields of the form $X=b(t)\partial/\partial x$ with $b$ being a solution of  (\ref{oscil}) are also infinitesimal symmetries of the equation.

Similarly, we can consider the nonlinear Ermakov-Pinney (EP) differential equation
\begin{equation}\label{EPk}
\ddot x+p(t) \, x=\frac {k}{x^3}, \quad x\ne0,\quad  k\in\mathbb{R}.
\end{equation}
One can show (see Section \ref{sec2}) that such differential equation is invariant under a 3-dimensional Lie algebra  of  Lie symmetries generated by vector fields of the form \eqref{EP-svfa} where $a$ satisfies Eq. \eqref{3rd2}. Let us mention that the differential equation \eqref{3rd2} is very related to the theory of higher order Adler-Gelfand-Dikii  differential operators \cite{A79,GD} and it plays a key role in the study of projective connections and $\mathfrak{gl}(n,\mathbb{R})$ current algebras\cite{M88, Gu02}.

The main objective of this paper is to  identify  families of second order ordinary differential equations that are invariant under a two-dimensional affine Lie subalgebra of the  Lie algebra associated to the EP equation, i.e., the Lie algebra generated by vector fields \eqref{EP-svfa}. One of the identified families of equations will be seen to include the EP equation \eqref{EPk} as a special case. This analysis is performed on Section \ref{sec3}.  The presence of the non-Abelian  two-dimensional symmetry Lie algebra, which is solvable, is sufficient for a second order ODE to be fully integrable  by quadratures. In a final section devoted to conclusions and remarks, we point out that the identified invariant  equations are not only of theoretical interest, but they are related to some recent models arising in population dynamics. This connection has been explored in more detail in a separate paper \cite{GuengoerTorres2019}.

\section{Preliminaries}\label{sec2}
In order to show our main motivation we will start with the derivation of  the Lie  algebra of infinitesimal point transformations of symmetry of  \eqref{EPk} using the prolongation algorithm for differential equations  (see for example \cite{Olver1993, Guengoer2019}).
Given the vector field  $X\in \mathfrak{X} (\mathbb{R}^2)$ with coordinate expression
\begin{equation}
X=\xi(t,x)\pd{}{t}+\eta (t,x)\pd{}{x}, \label{vfR2}
    \end{equation}
its second order prolongation $X^{(2)}$  is given by
\begin{equation}
 X^{(2)}= X +\eta
^{(1)}(t,x,\dot{x})\pd{}{\dot{x}}
 + \eta^{(2)}(t,x,\dot{x},\ddot{x})\pd{}{\ddot{x}},     \label{kpvfR2}
\end{equation}
where
$$\eta^{(1)}=D_t\eta-\dot{x}D_t\xi,  \quad \eta^{(2)}=D_t\eta^{(1)}-\ddot{x}D_t\xi.$$
Here  $D_t=d/dt$  is a symbol for the total derivative
\begin{equation}
\frac {d}{dt}=\pd{}{t} + \dot{x} \pd{}{x}+ \ddot{x}\pd{}{\dot{x}}+\cdots.\label{Ddt}
\end{equation}
More explicitly,  the coefficients of two first prolongations of (\ref{vfR2}) are:
$$
\eta^{(1)}=\eta_t+(\eta_x-\xi_t)\dot{x}-\xi_x\, \dot{x}^2,
$$
and
$$
D_t\eta^{(1)}=\eta_{tt}+(2\eta_{tx}-\xi_{tt})\dot{x}+(\eta_{xx}-2\xi_{xt})\dot{x}^2-\xi_{xx}\, \dot{x}^3-2\xi_x\, \dot{x}\ddot{x}+(\eta_x-\xi_t)\ddot{x},
$$
and consequently,
$$\begin{array}{rcl} \eta^{(2)}&=&D_t\eta^{(1)}-\ddot{x}D_t\xi=\eta_{tt}+(2\eta_{tx}-\xi_{tt})\dot{x}+(\eta_{xx}-2\xi_{xt})\dot{x}^2\\
&-&\xi_{xx}\, \dot{x}^3 +(\eta_x-2\xi_t)\ddot{x}-3\xi_x\, \dot{x}\ddot{x}.
\end{array}
$$

In the particular case of Ermakov-Pinney equation (\ref{EPk}) with $k\ne 0$,  the property characterizing the functions $\xi$ and $\eta$ such that  a vector field,
 (\ref{vfR2})   is a Lie symmetry of  such equation is given by
\begin{equation}\label{EPkeqsym}
\left(X^{(2)}(\ddot{x}+p(t)x-k\, x^{-3})\right)\Big|_{\ddot{x}+p(t)x-k\, x^{-3}=0}=0,\qquad  x>0,
\end{equation}
or more explicitly,
\begin{equation}\label{EPkeqsym2}
\left(\eta^{(2)}\right)\Big|_{\ddot{x}+p(t)x-k\, x^{-3}=0}+\dot p(t)x\xi+\eta\left(p(t)+3kx^{-4}\right)=0.
\end{equation}

The particularly interesting case is when the vector field is a projectable  vector field,
i.e. like in  (\ref{vfR2}) but with $\xi_x=0$,  because its flow is made of bundle map diffeomophisms $\phi_t:\mathbb{R}^2\to \mathbb{R}^2$.

The coefficients of the different powers of $\dot{x}$ in (\ref{EPkeqsym2}) must be zero, i.e. taking into account the corresponding form of $\eta^{(2)}$:
$$
\begin{array}{rl}(k\, x^{-3}-p(t)x)\!\!&\!\!\left(\eta_x-2\xi_t-3\xi_x\dot{x}\right)+\eta_{tt}+(2\eta_{xt}-\xi_{tt})\dot{x}
+(\eta_{xx}-2\xi_{xt})\dot{x}^2 -\xi_{xx}\dot{x}^3 \\&+\dot p(t)x\xi +\eta(p(t)+3kx^{-4})=0,
\end{array}
$$
and consequently
we find the following set of conditions:
\begin{equation}\label{conditions2}
\begin{array}{rl}
&\xi_{xx}=0,\\
&\eta_{xx}-2\xi_{tx}=0,\\
&2\eta_{tx}-\xi_{tt}+3(p(t)x-k\, x^{-3})\xi_x=0,\\
&(k\, x^{-3}-p(t)x)(\eta_x-2\xi_t)+\eta_{tt}+\dot p(t) x\xi+(p(t)+3kx^{-4})\eta=0.
\end{array}
\end{equation}

The two first equations lead to the following form for $\xi$ and $\eta$
$$
\xi(t,x)=d(t) x+b(t), \qquad \eta(t,x)= \dot d(t)x^2+c(t)\ x+e(t),
$$
and using these
expressions in the third equation of the preceding system we find
\begin{equation}\label{firstcondb}
2(2 \ddot d(t)x+\dot c(t))-(\ddot d(t)\, x+\ddot b(t))+3(p(t)x-k\, x^{-3}) d(t)=0.
\end{equation}
This condition implies, first,  that the function $d$ must be zero, because the coefficient of $x^{-3}$ is $k\,d(t)$,
and furthermore $2\dot c(t)-\ddot b(t)=0$, and then the expressions of the functions $\xi$ and $\eta$ are
$$
\xi(t,x)=b(t), \qquad \eta= c(t)x+e(t),
$$
which shows that $X$ is a projectable vector field.

Finally, the fourth equation reduces to
$$
(k\,x^{-3}-p(t)x)(c(t)-2\dot b(t))+\ddot c(t)x+\ddot e(t)+\dot p(t)x b(t)+(p(t)+3kx^{-4})(c(t)x+e(t))=0,
$$
and  for the coefficients of different powers of $x$ to be zero we obtain:
\begin{equation}\label{conditions3}
\begin{array}{rl}
&e(t)=0,\\
&2k\,(2c(t)-\, \dot b(t))=0,\\
&\ddot c(t)+2p(t)\, \dot b(t) +\dot p(t)\,  b(t)=0.
\end{array}
\end{equation}
 The second equation shows that
 \begin{equation}\label{cbdot}
 c(t)=\frac 12 \dot b(t),
 \end{equation}
 and a substitution in the third equation gives rise to
 \begin{equation}
 \frac 12 \dddot b(t)+2p(t)\, \dot b(t) +\dot p(t)\, \dot b(t)=0,
\label{eqb2}
\end{equation}
 i.e. $b$ is a solution of (\ref{3rd2}).

This means that  the symmetry vector fields we are looking for are of the form
 \begin{equation}
X_b(t,x)= b(t)\pd{}t+ \frac 12 \dot b(t)\,x\pd{}x,\label{XsymEP}
 \end{equation}
where $b(t) $ is a solution of (\ref{3rd2}).

The correspondence $a\mapsto X_a$ mapping each solution of (\ref{3rd2}) into  an infinitesimal symmetry of the nonlinear Ermakov-Pinney differential
equation (\ref{EPk})  is $\mathbb{R}$-linear, because $X_{a_1+\lambda\, a_2}=X_{a_1}+\lambda\, X_{a_2}$, for each real number $\lambda\in \mathbb{R}$.
Consequently, as (\ref{3rd2}) is a linear third order  differential equation, the set of vector fields  determined by solutions $a$ of the differential equation  (\ref{3rd2}) is a three-dimensional real linear space.

Note that if we consider Hill's equation (\ref{oscil}),
 it is possible to show that
if $u_1$ and $u_2$ are two linearly independent solutions of (\ref{oscil}), then the three functions $f_{ij}=u_i\,u_j$, $i\leq j=1,2$,  are  solutions of (\ref{3rd2}).

In fact, remark first that taking derivatives we obtain that
$$
  \dddot{u}_i+p(t)\,\dot u_i+\dot p(t)\,u_i = 0,  \quad i=1,2,
$$
and if we make use of these two equations, then the following third-order derivative
$$
 D_t^3(u_iu_j)=\dddot u_iu_j+3\ddot u_i\, \dot u_j+3\dot u_i\ddot u_j+u_i\dddot u_j\,,
$$
can be rewritten as follows
$$
 D_t^3(u_iu_j)=-(p(t)\,\dot u_i+\dot p(t)\,u_i)u_j-3p(t)\,u_i\dot u_j+3\dot u_i(-p(t)\,u_j)-
u_i(p(t)\,\dot u_j+\dot p(t)\,u_j),
$$
that after simplification becomes
$$
 D_t^3(u_iu_j)=-[2\dot p(t) \,u_iu_j+4p(t)(\dot u_iu_j+u_i\dot u_j)] \,.
$$
We have therefore obtained
$$
 D_t^3(u_iu_j)+4p(t)D(u_iu_j)+2\,\dot p(t)\,u_iu_j =0\,,
$$
what proves that the three functions $f_{ij}=u_i\,u_j$, $i\leq j=1,2$,  are  solutions of (\ref{3rd2}). Moreover, as the Wronskian of the three  functions $f_{ij}$ is
$$
W(u_1^2, u_1u_2, u_2^2) = 2 (u_1\,\dot u_2-u_2\,\dot u_1)^3 \,,
$$
we see that if $\{u_{1},u_{2}\}$ is a fundamental set of solutions of the second-order equation (\ref{oscil}), then
the functions $u_{1}^2$, $u_1 u_2$ and  $u_{2}^2$ are linearly independent
and they  span the  three-dimensional linear space of  solutions of (\ref{3rd2})
whose general solution can be written as a linear combination
\begin{equation}\label{sol-1}
  a(t)=A u_1^2+2Bu_1u_{2}+Cu_{2}^2,\qquad A,B,C\in \mathbb{R}.
\end{equation}

We can prove now that the set of vector fields of the form \eqref{EP-svfa} that are
Lie symmetries of the Ermakov-Pinney (EP) equation (\ref{EPk})
is   a Lie algebra: such
Lie symmetries of \eqref{EPk} close on the three-dimensional real Lie algebra $\Sl(2,\mathbb{R})$ spanned by the vector fields
\begin{equation}\label{EP-Lie}
  X_{ij}=f_{ij}\pd{} t+\frac{1}{2}\dot{f}_{ij}x\pd{} x,  \quad i\leq j=1,2.
\end{equation}

In fact,  the set of vector fields as in (\ref{EP-svfa}) is closed under commutator because
\begin{equation}
[X_{a_1},X_{a_2}]=X_{W(a_1,a_2)},\label{XaXa}
\end{equation}
where $W(a_1,a_2)$ denotes the Wronskian $W(a_1,a_2)=a_1\, \dot a_2-a_2\,\dot a_1$, and, moreover,  if $a_1$ and $a_2$ are solutions of (\ref{3rd2}), then the function
$w_{12}(t)=W(a_1(t),a_2(t))$ is a solution of  (\ref{3rd2}) too, because
$$(W(a_1,a_2))^{\dot{}}=\dot a_1\,\dot a_2+a_1\,\ddot a_2-\ddot a_1\,a_2-\dot a_2\,\dot a_1=a_1\,\ddot a_2-\ddot a_1\,a_2,
$$
and then,
$$
\dot w_{12}=(W(a_1,a_2))^{\dot{}}=a_1\,\ddot a_2-a_2\,\ddot a_1,
$$
and when taking derivatives in this  expression we get
$$
\ddot w_{12}=a_1\,\dddot a_2-a_2\,\dddot a_1+\dot a_1\,\ddot a_2- \dot a_2\,\ddot a_1,
$$
and
 therefore, if $a_1$ and $a_2$ are solutions  of (\ref{3rd2}),  a simple calculation shows that the preceding relation reduces to
 $$
 \ddot w_{12}=-4p(t) \, w_{12}+\dot a_1\,\ddot a_2- \dot a_2\,\ddot a_1,
 $$
 from where we see that
$$
\dddot w_{12}=-4p(t) \dot w_{12}-4\dot  p(t) \, w_{12}+\dot a_1(-4p(t)\,\dot a_2-2\dot p(t)\,a_2 )-\dot a_2(-4p(t)\,\dot a_1-2\dot p(t)\,a_1 ),
$$
and simplifying terms we arrive at
$$
\dddot w_{12}=-4p(t) \dot w_{12}-2\dot  p(t) \, w_{12}.
$$
We note that the same argument with more computational efforts can be used to show that  the Wronskian $w_{12}$ of any two independent solutions of  the general third order linear PDE
\begin{equation}\label{3rd-gen}
  \dddot{a}+c_2(t)\ddot{a}+c_1(t)\dot{a}+c_0(t)a=0,
\end{equation}
is also a solution
if and only if the coefficients satisfy $c_2=0$, $\dot{c}_1=2c_0$, (a formally self-adjoint equation).

Having in mind the mentioned property that for any pair of functionally independent solutions of  (\ref{oscil}), $u_1$  and $u_2$, the functions $u_1^2, u_1\, u_2$, and $ u_2^2$  form a  basis of the three-dimensional real linear space of solutions of (\ref{3rd2}),
we can consider as a basis of the three-dimensional real Lie algebra of infinitesimal symmetries of (\ref{3rd2}) the vector fields
$X_{u_1^2},  X_{u_1\,u_2},$ and $X_{u_2^2}$, and as
$$\begin{array}{rcl}
W(u_1^2, u_1\,u_2)&=&u_1^2\,W(u_1,u_2),\\ W(u_1^2,u_2^2)&=&2u_1\,u_2\,W(u_1\,u_2), \\   W(u_1\,u_2,u_2^2)&=&u_2^2\,W(u_1,u_2),
\end{array}
$$
where $W(u_1(t),u_2(t))$ is constant, and we obtain from (\ref{XaXa}) that
\begin{equation}\label{comm-sl2-2}
\begin{array}{rl}
&[X_{u_1^2},X_{u_1\,u_2}]=X_{W(u_1^2, u_1\,u_2)}=X_{u_1^2\,W(u_1, u_2)}, \\
&[X_{u_1^2},X_{u_2^2}]=X_{W(u_1^2, u_1\,u_2)}=2\,X_{u_1\, u_2\,W(u_1, u_2)},\\
& [X_{u_1\,u_2},X_{u_2^2}]=X_{W( u_1\,u_2,u_2^2)}=X_{u_2^2\,W(u_1, u_2)}.
\end{array}
\end{equation}

We can conclude from here that: {\it If $u_1$  and $u_2$ are two functionally independent solutions of (\ref{oscil}) such that $W(u_1, u_2)=1$,
then the vector fields $Y_1= X_{u_1^2}$,  $Y_2= X_{u_1\,u_2}$ and $Y_3= X_{u_2^2}$ generate a Lie algebra of vector fields of infinitesimal Lie symmetries of (\ref{3rd2}) isomorphic to $\mathfrak{sl}(2,\mathbb{R})$}, because they satisfy
the commutation relations
\begin{equation}\label{comm-sl2}
  [Y_1,Y_2]=Y_1,  \quad [Y_1,Y_3]=2Y_2,  \quad [Y_2,Y_3]=Y_3.
\end{equation}
This leads to the following result: The set of infinitesimal symmetries of (\ref{EPk})  is a three-dimensional real Lie algebra of vector fields like (\ref{EP-svfa}) where $a$ is
solution of (\ref{3rd2}).

It is also to be remarked that it has been proved in \cite{CarinenaLucasRanada2008a} that  the Ermakov-Pinney equation
$$
\ddot x=-\omega^2(t)x+\frac k{x^3},
$$ when written as a first-order system
$$\dot{x}=v,  \quad \dot{v}=-\omega^2(t)x+\frac k{x^3}$$
is a Lie system with associated Lie algebra  $\Sl(2,\mathbb{R})$, generated by the vector fields
$$X_1=x\, \pd{} v,\quad X_2=v\pd{} x+\frac{k}{x^3}\pd{} v,\quad X_3=\frac12 \left(x\pd{} x-v\pd{} v\right),
$$
which satisfy the following commutation relations
$$[X_1,X_2]= 2X_3,\qquad [X_1,X_3]=-X_1 ,\qquad [X_2,X_3]= X_2.
$$
As it happens for each Lie system, the flow of generators of its Vessiot Lie algebra,
the vector fields $X_1,X_2$ and $X_3$, transforms each  Lie system defined by them into another one of the same type.

Eq. \eqref{3rd2} admits the first integral
\begin{equation}\label{first-int}
  K=\frac{1}{4}(2a\,\ddot{a}-\dot{a}^2)+p(t) a^2,
\end{equation}
because multiplying the left hand side of (\ref{3rd2})  by $\frac 12 a$ we obtain
$$\frac 12 a\left(\dddot{a}+4p(t)\dot{a}+2\dot{p}(t) a\right)=\frac {d}{dt}\left(\frac{1}{4}(2a\,\ddot{a}-\dot{a}^2)+p(t)a^2\right)=0 .
$$

The value of $K$ for the general solution of (\ref{3rd2}) written in terms of two linearly independent solutions of (\ref{oscil})  as in \eqref{sol-1}  is specified as $K=(AC-B^2)w_{12}^2$, because
introducing the notation for the bilinear form $\<\cdot,\cdot>$
$$a=A u_1^2+2Bu_1u_{2}+Cu_{2}^2=(u_1,u_2)\begin{pmatrix}A&B\\B&C\end{pmatrix}\begin{pmatrix}u_1\\u_2\end{pmatrix}
\equiv\<\bu,\bu>,$$
where $\bu=u_1 \,\be_1+u_2\,\be_2$,  with $A=\<\be_1,\be_1>$, $B=\<\be_1,\be_2>$ and $C= \<\be_2,\be_2>$, then
$$\dot a=2 \<\dot \bu,\bu>,\qquad \ddot a= -2p(t)\<\bu,\bu>+2\<\dot\bu,\dot \bu>,
$$
and using the expression (\ref{first-int}) we find
$$K=\<\bu,\bu>(\<\dot\bu,\dot \bu>-p(t)\<\bu,\bu>-\<\dot\bu, \bu>^2)+p\<\bu,\bu>^2=\<\dot\bu,\dot \bu>\<\bu,\bu>-\<\dot\bu,\bu>^2.
$$

The right hand side of the preceding expression reminds that of the square of exterior product when
 $\<\cdot,\cdot>$  is the Euclidean product.  We can then define a skew-symmetric bilinear form  $F$
 either by this expression for the module  when the two vectors have positive orientation and the opposite  if the pair of vectors have the inverse orientation. This expression  $K=|F(\bu,\dot\bu)|^2$
shows that as for the exterior product
$$\|\bu\times \dot\bu\|=|W(u_1,u_2)| \ \|\be_1\times \be_2\|,$$
and any two skew-symmetric forms are proportional
$$K=|W(u_1,u_2)|^2\ |F(\be_1, \be_2)|^2,
$$
and as
$$|F(\be_1, \be_2)|^2=\<\be_1,\be_1>\<\be_2,\be_2>-\<\be_1,\be_2>^2=AC-B^2,$$
we find  from here the announced result.

We refer the interested readers to \cite{GuengoerTorres2017,CarinenaLucas2008} for solutions and Lie symmetry properties of EP
equation \eqref{EPk} and projective vector field equation  \eqref{3rd2}. Let us comment that Eq. \eqref{oscil} still has a  $\Sl(2,\mathbb{R})$  Lie algebra of Lie symmetry. In general, the third-order auxiliary equation \eqref{3rd2} crops up in symmetry analysis of second and higher order  linear ODEs with the property of being anti-self adjoint or
of maximal Lie symmetry and, as we have seen above, in second order nonlinear  ODEs whose solutions are expressed in terms
of \eqref{oscil} like EP and also its generalisations  \cite{Campoamor-Stursberg2016, GuengoerTorres2017},  and it is used in the  derivation of
 first integrals for time-dependent Hamiltonian systems \cite{StruckmeierRiedel2002}.

While it is possible to remove the coefficient $p$ from \eqref{EPk} by the change of variables $(t,x)\to (\tau,\xi)$ defined by
\begin{equation}\label{change-var}
  x=\xi(\tau)u_1,  \qquad \tau=(W(u_1,u_2))^{-1}\frac{u_2}{u_1},
\end{equation}
with  $u_1$ and $u_2$ particular solutions of (\ref{oscil}), we prefer to keep the potential $p$ to serve our purposes in the current context. Moreover, we can use  the orientation-preserving transformation
\begin{equation}\label{LF-tr}
  \bar{t}=\tau(t),  \quad \bar{a}(\bar{t})=\dot{\tau}(t)a(t),  \quad \dot{\tau}> 0,
\end{equation}
where $\tau$ satisfies the third-order  Kummer--Schwarz equation
\begin{equation}\label{Schwarz}
  \curl{\tau;t}=2\, p(t),
\end{equation}
with $\{\tau;t\}$ being the Schwarz derivative (see \cite{OT} for a short introduction), i.e.
\begin{equation}\label{Schwarzder}
\{\tau;t\}=\frac{\dddot{\tau}}{\dot{\tau}}-\frac{3}{2}\left(\frac{\ddot{\tau}}{\dot{\tau}}\right)^2.
\end{equation}
See e.g. \cite{CGL2013} and references therein, and \cite{deLucasSardon2013,CGLS2014,MPly} for related concepts and their physical applications.
Such transformation maps Eq.\eqref{3rd2} into its Laguerre-Forsyth canonical form $\bar{a}'''(\bar{t})=0$ \cite{neuman},  where the prime denotes derivative
 with respect to the new independent variable $\bar t$ (See for example \cite{homann,Olver1995}).  As remarked by Kummer \cite{Kummer} the solutions of \eqref{Schwarz} can be expressed as the quotient  of two linearly independent solutions of \eqref{oscil}.
This implies that transformation \eqref{LF-tr} can be written in the  form
\begin{equation}\label{trans2cano}
  \bar{t}=\tau(t)=\frac{\alpha u_1+\beta u_2}{\gamma u_1+\delta u_2},  \quad \bar{a}(\bar{t})=-\Delta\, W(u_1,u_2)(\gamma u_1+\delta u_2)^{-2}a(t),  \quad \Delta=\alpha \delta-\beta \gamma\ne 0.
\end{equation}
With the special choice $\alpha=0$, $\beta=1$, $\gamma=W(u_1,u_2)$, $\delta=0$ ($\Delta=-W(u_1,u_2)\ne 0$) and the relationship $x=\sqrt{a}$ between \eqref{first-int} and the equation $\ddot{x}+px=Kx^{-3}$, transformation \eqref{change-var} is recovered.

We can reobtain the general solution \eqref{sol-1} from \eqref{trans2cano}
$$a(t)=-\frac{1}{\Delta\, W(u_1,u_2)}(\gamma u_1+\delta u_2)^2 (c_1 + c_2 \tau+c_3 \tau^2)=A u_1^2+B u_1u_2+C u_2^2$$ after a redefinition of the arbitrary constants.

In particular, if we choose $p=0$ ($u_1(t)=t$, $u_2(t)=1$, $W=-1$) then we obtain the $\SL(2,\mathbb{R})$ subgroup of the  symmetry group of the canonical equation $\dddot{a}=0$
\begin{equation}
\bar{t}=\frac{\alpha t+\beta }{\gamma t+\delta },  \quad \Delta=1,\label{homot}
\end{equation}
together with $\bar{a}=(\gamma t+\delta)^{-2}a$.

The $\SL(2,\mathbb{R})$ symmetry group of the canonical EP equation $\ddot{x}=kx^{-3}$ is thus given by (\ref{homot})
with $  \bar{x}=(\gamma t+\delta)^{-1}x$.

In passing, we comment that Eq. \eqref{3rd-gen} can be reduced to the canonical form $\dddot{a}=0$ by a point transformation if and only if the following  singular invariant equation  relative to the general form-preserving transformation $\tau=\tau(t)$, $\bar{a}=\phi(t)a$ of \eqref{3rd-gen} is satisfied \cite{KrauseMichel1988,ndogmo}
\begin{equation}\label{invar}
  9\ddot{c}_2+18\dot{c}_2c_2-27\dot{c}_1+4c_2^3-18c_1c_2+54c_0=0.
\end{equation}
The special case $c_2=0$ is equivalent to the formal self-adjointness of the equation.

\section{Second order ODEs invariant under the two-dimensional affine algebra}\label{sec3}

We start  this section by looking for the second order  differential equations which admit as a Lie algebra of symmetry a Lie subalgebra of the Lie algebra of symmetries of the Ermakov-Pinney equation.
The only two-dimensional Lie subalgebra is isomorphic to that of the affine group of transformations of the real line. It is spanned by two vector fields $X_1$ and $X_2$ such that $[X_1,X_2]=X_1$.
Then, if $X_1$ and $X_2$ are vector fields of the form
$$X_{a_1}(t,x)=a_1(t)\pd{}t+\frac 12 \dot a_1(t)\,x\pd {}x,\qquad X_{a_2}(t,x)=a_2(t)\pd{}t+\frac 12 \dot a_2(t)\,x\pd {}x,
$$
where $a_1$ and $a_2$ are positive solutions of (\ref{3rd2}), then using the relation (\ref{XaXa})  we see that the functions $a_1$ and $a_2$ must be related by
$$W(a_1,a_2)=a_1,
$$
and therefore,
$$a_1\,\dot a_2-a_2\,\dot a_1=a_1,
$$
then, starting from a solution $a_1$ of (\ref{3rd2}) we obtain that $a_2$ must be a solution of the inhomogeneous linear differential equation
$$\dot  a_2=\frac{\dot a_1}{a_1}a_2+1.
$$
As $a_2=a_1$ is a solution of the associated linear homogeneous equation we should introduce the change of variable
$a_2=a_1\, s$, and the given equation becomes
$$a_1\,\dot s=1,
$$
which gives
$$
s(t)=\int^t\frac 1{a_1(\zeta)}\,d\zeta.
$$
Since $a_2=a_1\, s$ and
$$\frac d {dt}(a_1\, s)= \dot a_1\, s+a_1\, \dot s=  \dot a_1\, s+1,
$$
this proves that
$$X_{a_2}= s(t)\,X_{a_1}+\frac 12 x\,\pd{}x.
$$

We are now interested in  the most general class of second order ODEs involving functions expressed in terms of arbitrary solutions of \eqref{3rd2} for a given $p(t)$ and solvable by a pair of quadratures.
We start by realising the two-dimensional non-Abelian Lie algebra, generated by two vector fields $X_1$ and $X_2$ such that $X_1$ is  of the form \eqref{EP-svfa}, i.e.
is an infinitesimal point transformation  of symmetry of both (\ref{oscil}) and (\ref{EPk}),   and $X_2$ is a vector field  satisfying the  commutation relations
$[X_1,X_2]=X_1$.   Such Lie algebra is generated by
\begin{equation}\label{A2}
  X_1=a(t)\pd{}t+\frac{\dot{a}(t)}{2}x\pd{}x,  \quad X_2=s(t)X_1+\beta X_0, \quad X_0=x\pd{}x, \quad s(t)=\int^t \frac{d\zeta}{a(\zeta)},
\end{equation}
with  $\beta\ne 0$ and where $a$ was assumed to be  solution of \eqref{3rd2} for the  given $p(t)$  and $\beta$  a real number. This is so because given $X_1$ of the above mentioned  form, then we can write $X_2$ as a linear combination
 of the form
 $$X_2=c(t)\, X_1+b(t)\, x\pd{}x=c(t)\, X_1+b(t)\,X_0,
$$
and then, as $[X_1,X_0]= 0$,
$$
[X_1,X_2]=[X_1,c(t)\,X_1]+\left[X_1,b(t)\, x\pd{}x\right]= X_1(c)\, X_1+a(t)\, \dot b(t)\, x\pd{}x ,
$$
and therefore, in order to have $[X_1,X_2]=X_1$, the functions $b$ and $c$ must satisfy
$$
a\dot c=1, \qquad \dot b=0,
$$
from where we obtain that $b(t)$ must be constant, $b(t)=\beta$, and  $c(t)$ must be given by $s(t)$ as indicated by (\ref{A2}). The constant $\beta$ must be different from zero,
otherwise $X_2$ and $X_1$ would be proportional in each point

The 2-dimensional  Lie algebra spanned by $X_1$ and $X_2$ (recall that we assumed $\beta\ne 0$)  is isomorphic to the Lie algebra of the affine group  in the real line.
They define a  transitive  action of this Lie algebra  on the  plane $(t,x)$  (no nontrivial ordinary invariants exist). Recall that if $a$ is not constant, only in the  particular case  $\beta=1/2$
the vector field $X_2$ is of the family of vector fields (\ref{EP-svfa}), in other words we only have a rank-two realisation of the algebra within the class of vector fields \eqref{EP-svfa}.

Our aim is to construct the general second order ODE invariant under the realisation \eqref{A2} of the  two-dimensional  affine algebra $\aff(1,\mathbb{R})$. This is a standard procedure and requires finding  invariants for the second prolongation $\pr{2}\aff(1,\mathbb{R})$ by solving a pair of first order linear PDEs  by the method of characteristics.

\subsection{Invariant equation in the case $\beta=1/2$.}

We should look for the most general second order ODE invariant under the realisation \eqref{A2} with $\beta=1/2$. We start by looking for
the invariant functions  for the second prolongation of  $X_1$,  $X_1^{(2)}$, given by
$$
X_1^{(2)}= a(t)\pd{}t+\frac{\dot{a}(t)}{2}x\pd{}x+\frac{1}{2}\left(\ddot{a}x-\dot{a}\dot{x}\right)\pd {}{\dot x}+\frac{1}{2}\left (\dddot{a}x-3\dot{a}\ddot{x}\right)\pd {}{\ddot x}.
$$
 The can be are computed as characteristic solutions of the partial differential equation $X_1^{(2)}H=0$
 and we find as solution a function $ H(J_1,J_2,J_3)$ where
\begin{equation}\label{inv-1}
  J_1=\frac{x}{\sqrt{a}}, \quad J_2=\sqrt{a}\left(\dot{x}-\frac{\dot{a}}{2a}x\right)=a\dot{J}_1, \quad
  J_3=a^{3/2}(\ddot{x}+px).
\end{equation}
In order to impose  $X_2$ invariance we first remark that as
$$X_2^{(2)}(J_1)=\frac{1}{2}J_1,  \quad X_2^{(2)}(J_2)=-\frac{1}{2}J_2$$ and
$$X_2^{(2)}(J_3)=\frac{1}{2}a^{3/2}(-3\ddot{x}+px)+\frac{a^{-1/2}}{2}(2a\ddot{a}-\dot{a}^2)x,$$
that using the first integral \eqref{first-int} we can rewrite as
$$X_2^{(2)}(J_3)=-\frac{3}{2}J_3+2KJ_1,$$
the differential invariants of order $\leq 2$ of  the algebra $\aff(1,\mathbb{R})$ are found by solving the PDE
\begin{equation}\label{PDE}
  \frac{1}{2}J_1\frac{\partial H}{\partial J_1}-\frac{1}{2}J_2\frac{\partial H}{\partial J_2}+\left(-\frac{3}{2}J_3+2KJ_1\right)\frac{\partial H}{\partial J_3}=0.
\end{equation}
Then we consider the associated system
$$\frac{dJ_1}{J_1}=-\frac{dJ_2}{J_2}=\frac{dJ_3}{-3J_3+4KJ_1}
$$
and  from the first fraction with the second or with the third one we find the invariant  functions
\begin{equation}\label{inv-2}
  I=J_1J_2=x\dot{x}-\frac{\dot{a}}{2a}x^2,  \quad J=J_1^3 J_3-K J_1^4=x^3(\ddot{x}+px)-\frac{K}{a^2}x^4,
\end{equation}
such that the general solution of  $X_1^{(2)}H=X_2^{(2)}H=0$ is an arbitrary function of $I$ and $J$.

The  invariant second order ODE are therefore of the form
\begin{equation}\label{inv-eq-0}
  x^3(\ddot{x}+px)=\frac{K}{a^2}x^4+G(I),
\end{equation}
with $G$ an arbitrary smooth function, or written in a different way,
\begin{equation}\label{inv-eq-0-1}
  \ddot{x}+\left(p(t)-Ka(t)^{-2}\right)x=x^{-3}G(I).
\end{equation}
  So for a given function $p$ we can produce a class of  ODEs integrable by quadratures. The first integral condition \eqref{first-int} gives us
\begin{equation}\label{inv-eq-0-2}
  \ddot{x}-\frac{1}{4a^2}(2a\ddot{a}-\dot{a}^2)x=x^{-3}G\left(x\dot{x}-\frac{\dot{a}}{2a}x^2\right),
\end{equation}
which  actually depends on $p$ in a disguise form. The above equation can be written explicitly  with  $\nu=\dot{a}/a$ as
\begin{equation}\label{inv-eq}
  \ddot{x}-\frac{1}{4}(2\dot{\nu}+\nu^2)x=x^{-3}G\left(x\dot{x}-\frac{1}{2}\nu x^2\right).
\end{equation}
It is straightforward to see that Eq. \eqref{inv-eq} allows the following invariant (particular solutions)
\begin{gather}\label{part-sol}
  x(t)=C_0\sqrt{a(t)}, \qquad G(0)=0, \\
  x(t)=C_0\sqrt{s(t)a(t)},  \qquad C_0+4G\left(\frac{ C_0^2}{2}\right)=0.
\end{gather}
It is useful to find an equivalent form of \eqref{inv-eq} under the transformation $x=z^{1/k}$, which takes \eqref{inv-eq-0} to
\begin{equation}\label{inv-eq-2}
  \ddot{z}-\frac{k}{4}\left(2\dot{\nu}+\nu ^2\right)z=\frac{k-1}{k}\frac{\dot{z}^2}{z}+kz^{(k-4)/k}G(I),
\end{equation}
where
$$I=\frac{1}{k}z^{(2-k)/k}\left(\dot{z}-\frac{k}{2}\nu z\right), \quad \nu=\frac{\dot{a}}{a}.$$
Of course, the constant $k$ multiplying  the arbitrary function $G$ can be absorbed into $G$.  The symmetry algebra is
\begin{equation}\label{A2-k}
  X_1=a(t)\pd{}t+\frac{k\dot{a}(t)}{2}z\pd{}z,  \quad X_2=s(t)X_1+\frac{k}{2}z\pd{}z=a(t)s(t)\pd{}t+\frac{k}{2}(1+\dot{a}s)z\pd{}z.
\end{equation}
It is more convenient to put $k=4/(1-n)$ for some real $n\ne 1$ for which \eqref{inv-eq-2} takes the form
\begin{equation}\label{inv-eq-0-n}
\ddot{z}+\frac{1}{n-1}(2\dot{\nu}+\nu^2)z=\frac{n+3}{4}\frac{\dot{z}^2}{z}+z^{n}G(I),
\end{equation}
where
$$I=\frac{1-n}{4}z^{-(n+1)/2}\left(\dot{z}-\frac{2\nu}{1-n} z\right), \quad \nu=\frac{\dot{a}}{a}.$$

Now we will examine some particular cases. For $p=-\lambda^2/4$ we have the possibilities $a=1$ ($\nu=0$, $K=-\lambda^2/4$), and $a=e^{\pm \lambda t}$, $\lambda\ne 0$ ($\nu=\pm \lambda$, $K=0$) and the corresponding invariant equations have the form
\begin{gather}\label{ex-1}
  \ddot{x}=x^{-3}G(x\dot{x}), \\
  \ddot{x}-\frac{\lambda^2}{4}x=x^{-3}G\left(x\dot{x}\pm\frac{\lambda}{2}x^2\right).
\end{gather}

The corresponding symmetry vector fields are
\begin{equation}\label{X1-X2-1}
  \begin{gathered}
     X_1=\pd{}t,  \quad X_2=t\pd{}t+\frac{x}{2}\pd{}x, \\
     X_1=\exp[\pm \lambda t]\left(\pd{}t\pm \frac{\lambda}{2}x\pd{}x\right),  \quad X_2=\pm \frac{1}{\lambda}\pd{}t.
   \end{gathered}
\end{equation}

For $p=\lambda^2/4$, we have either $a=\cos (\lambda t)$ (and $\nu=-\lambda \tan (\lambda t)$, $K=-\lambda^2/4$) or $a=\sin (\lambda t)$ (and then $ \nu=\lambda \cot (\lambda t)$, $K=-\lambda^2/4$). The value of $K$ is determined either by direct computation from \eqref{first-int} or by making use of the relation $K=(AC-B^2)W^2$. For example, comparing the relation
$$a=\cos (\lambda t)=\cos^2 \frac{\lambda t}{2}- \cos^2 \frac{\lambda t}{2}=u_1^2-u_2^2$$ with \eqref{sol-1}
 implies $A=-C=1$, $B=0$ and with $W(u_1,u_2)=\lambda/2$ we find $K=-\lambda^2/4$. For $a=\sin (\lambda t)$, we have $A=C=0$, $B=1$.
The corresponding equation and symmetries for $a=\cos (\lambda t)$ are
\begin{equation}\label{ex-2}
 \ddot{x}+\frac{\lambda^2}{4}(1+\sec^2 (\lambda t)) x=x^{-3}G\left(x \dot{x}+\frac{\lambda}{2}\tan (\lambda t) \, x^2\right),
\end{equation}
\begin{equation}\label{X1-X2-2}
\begin{gathered}
  X_1=\cos (\lambda t)\pd{}t-\frac{\lambda}{2}\sin (\lambda t) x\pd{}x, \\
   X_2=\frac{2}{\lambda} \tanh^{-1}\left(\tan \frac{(\lambda t)}{2}\right)\cos (\lambda t)\pd{}t+\frac{1}{2}\left(1-2\tanh^{-1}\left(\tan \frac{\lambda t}{2}\right)\sin (\lambda t)\right)x\pd{}x.
\end{gathered}
\end{equation}

\subsection{Reduction to quadrature and solutions}
We can introduce the new coordinates $(r,s)$ adapted to the vector field $X_1$, i.e. such that  $X_1r=0, X_1s=1$, which are therefore given by
\begin{equation}\label{canon}
  r=\frac{x}{\sqrt{a}},   \quad s=\int^t \frac{d\zeta}{a(\zeta)},
\end{equation}
so that $X_1= \partial /\partial s$. Then as $X_2r=r/2$  and $X_2s=s$,  the rank-two affine algebra is transformed to the one generated by the vector fields $$X_1=\pd{}s, \quad X_2=s\pd{}s+\frac{r}{2}\gen r.$$

If we note the relations
$$r\frac{dr}{ds}=x\left(\dot{x}-\frac{\dot{a}}{2a}x\right),  \qquad \frac{d^2r}{ds^2}=a^{3/2}\left(\ddot{x}-\frac{1}{4a^2}(2a\ddot{a}-\dot{a}^2)x\right),$$
the canonical form of invariant equation \eqref{inv-eq-0-1} or \eqref{inv-eq-0-2}  has the  form of the generalised Ermakov-Pinney equation
\begin{equation}\label{trans-eq}
  \frac{d^2r}{ds^2}=r^{-3}G\left(r\frac{dr}{ds}\right).
\end{equation}
The equivalent form \eqref{inv-eq-0-n} is reduced to the  canonical form
\begin{equation}\label{trans-eq-2}
  r''(s)=\frac{n+3}{4}\frac{r'^{2}}{r}+r^{n}G(\omega), \quad  \omega=\frac{(1-n)}{4}r^{-(n+1)/2}r'
\end{equation}
by means of the coordinate transformation
\begin{equation}\label{canon-2}
  r=a^{2/(n-1)}z,   \quad s=\int^t \frac{d\zeta}{a(\zeta)}.
\end{equation}
Eq. \eqref{trans-eq-2} is invariant under the Lie algebra spanned by the vector fields $$X_1= \pd{} s, \quad X_2=s\pd{} s+\frac{2}{1-n}r\pd{} r.$$
When $G$ is restricted to a constant, say $G=4G_0/(1-n)$, with $G_0$ a constant, it is known as a special case of second order Kummer-Schwarz equation (see Eq. \eqref{stand}), which has a general solution formula so that  solution $z$ of \eqref{inv-eq-0-n} is given by
\begin{equation}\label{sol-z}
  z(t)=a^{2/(1-n)}r(s)=(Aa+2Bas+Cas^2)^{2/(1-n)},  \quad AC-B^2=G_0.
\end{equation}

The structure of the canonical equation \eqref{trans-eq} or \eqref{trans-eq-2} for $n=-3$ suggests the special choice $G=\text{const.}$ which reduces to the canonical form of the standard Ermakov-Pinney equation, namely
\begin{equation}\label{EP-cano}
  r''(s)=G_0r^{-3}.
\end{equation}
In this case the affine symmetry algebra \eqref{A2} with $\beta\ne 0$ extends to an $\Sl(2,\mathbb{R})$ algebra isomorphic to the second type in Lie's classification list. The additional symmetry vector field is given by
$$X_3=s^2\gen s+sr\gen r.$$
 We already know that Eq. \eqref{EP-cano} admits a general solution formula given by
\begin{equation}\label{sol-EP-cano}
  r(s)=(A+2Bs+Cs^2)^{1/2},  \quad AC-B^2=G_0.
\end{equation}
From this fact we immediately see that the following equation
\begin{equation}\label{inv-eq-sl2}
  \ddot{x}+[p(t)-Ka(t)^{-2}]x=G_0x^{-3}
\end{equation}
admits a $\Sl(2,\mathbb{R})$ symmetry algebra spanned by the vector fields \eqref{A2} and an additional one
\begin{equation}\label{X3}
X_3=a(t)s(t)^2\gen t+\frac{1}{2}(\dot{a}(t)s(t)^2+2s(t))x\gen x.
\end{equation}
We note that the  realisation of the $\Sl(2,\mathbb{R})$ Lie algebra is generated by
\begin{equation}\label{sl2-vf}
  X_1=a\gen t+\frac{\dot{a}}{2}x\gen x,  \quad X_2=sX_1+\frac{1}{2}X_0,  \quad X_3= s^2X_1+sX_0,
\end{equation}
where  the vector field      $X_0$ is $X_0=x\partial/\partial  x$,
with commutation relations
\begin{equation}\label{comm-sl2-3}
  [X_1,X_2]=X_1,  \quad [X_1,X_3]=2X_2,  \quad [X_2,X_3]=X_3,
\end{equation}
which can be derived from the following commutation relations
\begin{gather}\label{comm-rel}
  [X_1,X_0]=0,  \quad [X_1,sX_0]=X_0,  \quad [X_1,sX_1]=X_1,\\
  [X_1,s^2X_1]=2sX_1,  \quad [sX_1,sX_0]=sX_0,  \quad [sX_1,s^2X_1]=s^2X_1,
\end{gather}
from where we see that
 the symmetry vector fields $X_1,X_2,X_3$ of \eqref{sl2-vf} satisfy the commutation relations (\ref{comm-sl2-3}) characteristics of $\Sl(2,\mathbb{R})$  Lie algebra.
 The general solution of \eqref{inv-eq-sl2} is now given by the formula
\begin{equation}\label{sol-2}
  x(t)=\sqrt{a(t)}r(s(t))=\sqrt{Aa+2Bas+Cas^2}, \quad AC-B^2=G_0.
\end{equation}
This solution is somewhat surprising because as long as $K$ is a non-vanishing constant we obtain the general solution of the Ermakov-Pinney  equation with a considerably modified potential $\tilde{p}(t)=p(t)-Ka(t)^{-2}$,  and only when $K=0$ it coincides with the usual Ermakov-Pinney solution.

As an example we consider a case where $K=-\lambda^2/4\ne 0$, $p=\lambda^2/4$, $a=\cos (\lambda t)$:
\begin{equation}\label{ex-sl2}
  \ddot{x}+\frac{\lambda^2}{4}(1+\sec^2 (\lambda t))x=G_0x^{-3}.
\end{equation}
The general solution of \eqref{ex-sl2}, despite being too complicated, is  given exactly by  the formula \eqref{sol-2} with $s(t)$ being
$$s(t)=\frac{1}{\lambda}\log\left[\frac{1+\tan \frac{\lambda t}{2}}{1-\tan \frac{\lambda t}{2}}\right].$$
The choice  $p=1$, $a=1+\alpha \cos (2t)$, $|\alpha|<1$ ($K=1-\alpha^2$) leads to the EP equation
\begin{equation}\label{ex-sl2-2}
  \ddot{x}+\left(1+\frac{\alpha^2-1}{(1+\alpha \cos (2t))^2}\right)x=G_0x^{-3},  \qquad |\alpha|<1.
\end{equation}
The $\pi$-periodic general solution of \eqref{ex-sl2-2} is given by \eqref{sol-2}  with $s(t)$ being
$$s(t)=\frac{1}{\sqrt{1-\alpha^2}}\arctan\left(\frac{1-\alpha}{1+\alpha}\tan t\right).$$
The linear version of \eqref{ex-sl2-2} with $G_0=0$ belongs to a one-parameter family of Hill's equations with coefficients periodic of period $\pi$ (also a subclass of the so-called four-parameter Ince equations \cite{Athorne1990}).

On the other hand, the special choice $G(I)=4G_0/(1-n)$ in \eqref{inv-eq-0-n}  produces the following important form of a $\Sl(2,\mathbb{R})$-invariant equation that frequently arises in many applications
\begin{equation}\label{inv-eq-other}
  \ddot{z}+\frac{4}{1-n}\left(p-Ka^{-2}\right)z=\frac{n+3}{4}\frac{\dot{z}^2}{z}+\frac{4G_0}{1-n}z^n.
\end{equation}
A basis of the symmetry algebra is given by
    \begin{equation}\label{comm-sl2-KS}
    X_1=a\gen t+\frac{k\dot{a}}{2}x\gen x,  \quad X_2=sX_1+\frac{k}{2}X_0,  \quad X_3= s^2X_1+ksX_0,  \quad X_0=x\gen x,
 \end{equation}
  where $k=4/(1-n)$.
The general solution of \eqref{inv-eq-other} is given by (see solution \eqref{sol-z})
\begin{equation}\label{sol-3}
  z(t)=(Aa+2Bas+Cas^2)^{2/(1-n)},  \quad AC-B^2=G_0.
\end{equation}
This equation can be regarded as a generalisation of the second order Kummer-Schwarz (2KS) equation provided that $K\ne 0$.

The following dissipative form of \eqref{inv-eq-other} for $K=0$ can also be of some interest
\begin{equation}\label{diss-KS}
  \ddot{w}+r(t)\dot{w}+\frac{4p(t)}{1-n}w=\sigma\frac{\dot{w}^2}{w}+\frac{4q}{1-n}\exp\left[-2\int ^t r(\zeta)d\zeta\right]w^n, \ n\ne 1,\ q\in\mathbb{R}, \ \sigma=\frac{n+3}{4}.
\end{equation}
We call \eqref{diss-KS}  dissipative second order Kummer-Schwarz (d2KS) equation.
The linear transformation
\begin{equation}\label{trans}
  w(t)=\phi(t)z(t),  \quad \phi(t)=\exp\left[\frac{1}{2(\sigma-1)}\int^t r(\zeta)d\zeta\right],  \quad 2(\sigma-1)=\frac{n-1}{2},
\end{equation}
transforms \eqref{diss-KS} into
\begin{equation}\label{stand}
  \ddot{z}+\frac{4}{1-n}I(t)z=\sigma\frac{\dot{z}^2}{z}+\frac{4q}{1-n}z^n,
\end{equation}
where $$I(t)=p-\frac{1}{4}(r^2+2\dot{r}).$$
We already know that Eq. \eqref{stand} has the general solution
\begin{equation}\label{super}
  z=(Au_1^2+2Bu_1u_2+Cu_2^2)^{2/(1-n)},  \quad (AC-B^2)W^2(u_1,u_2)=q,
\end{equation}
where $u_1, u_2$ are two linearly independent solutions of the equation
\begin{equation}\label{base}
  \ddot{z}+I(t)z=\ddot{z}+\left(p-\frac{1}{4}(r^2+2\dot{r})\right)z=0.
\end{equation}
The general solution of  \eqref{diss-KS} is given by
\begin{equation}\label{gen-sol-dissip}
 w(t)=\exp\left[\frac{2}{n-1}\int^t r(\zeta)d\zeta\right](Au_1^2+2Bu_1u_2+Cu_2^2)^{2/(1-n)},  \quad (AC-B^2)W^2(u_1,u_2)=q.
\end{equation}

The d2KS equation \eqref{diss-KS} is invariant under   the real Lie algebra of vector fields
$$X_a=a\gen t+\frac{2}{1-n}(\dot{a}-ar)w\gen w,$$
where the function $a$ is in the real  linear space spanned   by the functions $u_1^2,u_1u_2,u_2^2$,  where  $u_1, u_2$ are solutions of
$$\ddot{w}+[p-\frac{1}{4}(r^2+2\dot{r})]w=0.$$
The commutation relations between the three components of the algebra  satisfy those of the $\Sl(2,\mathbb{R})$ algebra in \eqref{comm-sl2}.

We note that a Lagrangian $L$ of the 2KS equation \eqref{stand} is provided by
\begin{equation}\label{Lag}
  L(t,z,\dot{z})=\left(\frac{1-n}{4}\right)^2z^{-(n+3)/2}\dot{z}^2-I(t)z^{(1-n)/2}-qz^{(n-1)/2}.
\end{equation}

\subsection{Reduction to quadratures of Eq.   \eqref{trans-eq-2} }
We now turn to perform reduction   to quadratures of the differential equation \eqref{trans-eq-2}. To this end, we let $R=dr/ds$ and exchange the roles of $(r,s)$. This gives the first order equation
$$\frac{dR}{dr}=\frac{n+3}{4}\frac{R}{r}+\frac{r^n}{R}G(\omega), \quad \omega=\frac{1-n}{4}r^{-(n+1)/2}R.$$
 Invariance of this equation under the dilational symmetry generated by the vector field  $ r\partial/\partial r+\frac{(n+1)}{2}R\,\partial/\partial R $ implies
reduction  to the  separable form
\begin{equation}\label{sep}
  \frac{d\omega}{d\xi}=\frac{1-n}{4}\omega+\frac{(1-n)^2}{16\omega}G(\omega),
\end{equation}
which is achieved by changing coordinates to $(\omega,\xi=\ln r)$  and $r, s$ defined by \eqref{canon-2}. Once a solution $\omega=\Phi(\xi,C_1)$ to \eqref{sep} has been found,
the general solution is obtained by integrating another separable first order ODE
$$\frac{dr}{ds}=R=\frac4{(1-n)}r^{(n+1)/2}\Phi(\ln r,C_1).$$

More conveniently, one can use the change of coordinates $\bar{s}=r^{(1-n)/2}$, $\bar{r}=s+r^{(1-n)/2}$ to transform \eqref{trans-eq-2} into
\begin{equation}\label{canon-3}
  \bar{s}\frac{d^2 \bar{r}}{d\bar{s}^2}=\widehat{G}\left(\frac{d\bar{r}}{d\bar{s}}\right)
\end{equation}
with symmetry Lie algebra generated by $\langle \partial/\partial {\bar{r}}, \bar{s}\partial/\partial {\bar{s}}+\bar{r}\partial/\partial {\bar{r}}\rangle$ and a new arbitrary function $\widehat{G}$.
Integration of \eqref{canon-3} is straightforward.

\subsection{Linearizable subclasses by Lie's test}
In this subsection, we reconsider the canonical equation \eqref{trans-eq} for $r(s)$
\begin{equation}\label{cano-ode}
  r''=f(r,r')=r^{-3}G(rr')=r^{-3}G(I)
\end{equation}
and apply the Lie's test for a second order ODE in normal form $r''=f(s,r,p)$, $p=r'$, which determines the necessary and sufficient
conditions for transformability to a linear equation by a point transformation. Such conditions are expressed by the vanishing of the following
fourth order relative invariants \cite{MilsonValiquette2015}
\begin{equation}\label{abs-inv}
  \mathbb{I}_1=f_{pppp}=0,  \quad \mathbb{I}_2=\widehat{D}_s^2f_{pp}-4\widehat{D}_sf_{rp}-f_p \widehat{D}_s f_{pp}+6f_{rr}-3f_rf_{pp}+4f_pf_{rp}=0,
\end{equation}
where $\widehat{D}_s=\partial/\partial s+p\partial/\partial r+f\partial/\partial p$. The first condition requires that $G$  must be a cubic polynomial  of $I$, $G(I)=G_0 I^3+G_1 I^2+G_2I+G_3$. The second condition restricts the coefficients in two possible forms
\begin{gather}\label{linear-eq-conds}
  G_2=G_3=0, \\
  G_0=\frac{G_2}{27G_3^3}(G_2^2-18G_3),  \quad G_1=\frac{G_2^2-5G_3}{3G_3},  \quad G_3\ne 0.
\end{gather}
The first choice gives the equation $r''=G_0r'^3+G_1s^{-1}r'^2$, which is equivalent to the linear equation $s''(r)+G_1r^{-1}s'(r)+G_0=0$ by an exchange of the coordinates $s\leftrightarrow r$.

The other possibility gives the linearizable equation
\begin{equation}\label{linear-eq}
  r''(s)=\frac{G_2}{6G_3}\left(\frac{G_2^2}{18G_3}-1\right)\frac{r'^3}{r^3}+3\left(1-\frac{G_2^2}{18G_3}\right)
  \frac{r'^2}{r}+G_2rr'-2G_3r^3.
\end{equation}
Reverting $(s,r)$ back to $(t,x)$ gives us a more general form of a linearizable second order ODE.

The special choice $G_3=G_2^2/18$, $G_2=-3\ell$ of the coefficients singles out a well-known second member of the Riccati chain (the modified Emden equation) \cite{GrundlandLevi1999, CarinenaGuhaRanada2009}
\begin{equation}\label{2nd-Riccati}
  r''+3\ell rr'+\ell^2 r^3=0,
\end{equation}
which is generated by the second iteration of the Riccati operator $\mathbb{D}=D_s+\ell r$:
\begin{equation}\label{Riccati-chain}
  \mathbb{D}^2r=(D_s+\ell r)(D_s+\ell r)r=0.
\end{equation}
Eq. \eqref{Riccati-chain} is also recognised as a spacial case of the second order Riccati equation in the sense of Vessiot and Wallenberg \cite{Davis1962}. This $\Sl(3,\mathbb{R})$ invariant equation can also be obtained from \eqref{inv-eq-2} by choosing $k=-2$, $a=1$ ($s(t)=t$) and $G(I)=3I^2-3\ell I+\ell^2/2$.
By scaling  $r \to \ell r$ we can put $\ell=1$.

Just like the ordinary first order Riccati equation, the Hopf--Cole transformation $r=\rho'/\rho$ linearizes \eqref{2nd-Riccati} to the third order linear equation $\rho'''=0$.
Moreover, a point transformation linearizing \eqref{2nd-Riccati} to $R''(S)=0$ is provided by (see Example 5.5 of \cite{Guengoer2019})
\begin{equation}\label{lineariz-tr}
 S=s-\frac{1}{r},  \quad R=\frac{s^2}{2}-\frac{s}{r}.
\end{equation}
We comment that though the more general form
\begin{equation}\label{more-gen}
  r''+arr'+br^3=0
\end{equation}
does not pass Lie test unless $b=a^2/9$, it was shown to be linearizable to
$$\frac{d^2 \bar{r}}{d\bar{s}^2}+a\frac{d\bar{r}}{d\bar{s}}+2b\bar{r}=0$$
by the nonlocal transformation $\bar{s}={\displaystyle\int ^s}r(\zeta)\, d\zeta$, $\bar{r}=r^2$ \cite{Abraham-Shrauner1993}.

Finally, we mention that it was shown in \cite{CarinenaRanadaSantander2005}  using an ansatz that a special case of second-order Riccati equation, in particular \eqref{2nd-Riccati} with $\ell=1$, admits the (non-natural) Lagrangian
\begin{equation}\label{Lag-2}
  L=\frac{1}{r'+ r^2}.
\end{equation}
We can recover $L$ by transforming the Lagrangian  $L_0=R'^2$ of the free particle equation $R''=0$ by \eqref{lineariz-tr}. The transformed Lagrangian $\bar{L}$ is obtained as
$$\bar{L}=\frac{[s(r'+r^2)-r]^2}{(r'+r^2)^2}D_sS=
\frac{1}{r'+r^2}+s^2+\frac{s^2r'}{r^2}-\frac{2s}{r}=L+D_s\left[\frac{s^3}{3}-\frac{s^2}{r}\right].$$
Remark that as  the  Lagrangians $\bar{L}$ and $L$ differ by a total derivative they give rise to  the same Euler-Lagrange equation \eqref{2nd-Riccati}. In other words, $L$ and  $\bar{L}$ are gauge equivalent Lagrangians \cite{CarinenaIbort}.

\section{Conclusions and outlook}\label{sec4}
In this paper we have analysed the invariance of second order ODEs under a  2-dimensional affine Lie algebras realised by vector fields \eqref{A2} as extensions of the EP-symmetry vector field
 \eqref{EP-svfa}.
 By construction, these type of equations can be integrated by Lie's standard reduction procedure. It is also possible to give some particular (invariant) solutions.
 In the rank two case, for a constant choice
 of the arbitrary function $G$ appearing in the ODE, we have produced  an equation of EP type (see \eqref{EPk}) but  with potential $p(t)$ replaced by $p(t)-Ka^{-2}(t)$, $K$ being some constant fixed by choice of $a$.
  The general solution formula for  (\ref{EPk}) remains unchanged. We have introduced a dissipative
   (damped) version of EP equation and presented its general solution (nonlinear superposition). Linearisable subclasses of the canonical ODEs are obtained by Lie's test.

As a final remark, let us mention that the presented study is not merely academic, for some equations treated here arise in different applications. For example, in the recent paper \cite{NucciSanchini2015}, the authors investigated solutions and first integrals of a second order ODE falling within
the class \eqref{diss-KS}, based on symmetry approach. This ODE is obtained from elimination of a dynamical system modeling the total population of Easter island  \cite{BasenerRoss2004}. Solutions can be readily recovered  from our general results. A separate article \cite{GuengoerTorres2019} has recently been devoted to study integrability properties of a variable coefficient variant of the above-mentioned model by using results of the present work.

\bibliographystyle{unsrt}

\subsection*{Acknowledgments}
This work was initiated while one of the authors (F. G.) was visiting the Department of Theoretical Physics, Zaragoza University. This author is much indebted to J.~F. Cari{\~n}ena for invaluable discussions on Lie and Quasi-Lie systems and also to  members of the Department  for the warm hospitality.

\end{document}